\documentclass[reqno,12pt]{amsart}
\usepackage{indentfirst}
\usepackage{amsmath}
\usepackage[T1]{fontenc}
\usepackage{lmodern}
\usepackage{amssymb}			
\usepackage{amsthm}
\usepackage{amsfonts}
\usepackage{enumerate}
\usepackage{enumitem}
\usepackage{setspace}
\usepackage{mathtools}
\usepackage[body={7in,9.5in},top=0.5in, left=0.8in ]{geometry}
\newcommand{\bn}{\mathbb{N}}

\newcommand{\br}{\mathbb{R}}

\newcommand{\ov}[1]{\overline{#1}}
\newcommand{\wh}[1]{\widehat{#1}}

\newcommand{\pref}[1]{\textup{(\ref{#1})}}

\newtheorem{thrm}{Theorem}
\newtheorem{lm}[thrm]{Lemma}

\author[P.Ma\'ckowiak]{Piotr Ma\'ckowiak}
\address{P.Ma\'ckowiak, Faculty of Mathematics and Computer Science\\
Department of Nonlinear Analysis and Applied Topology\\
Adam Mickiewicz University in Pozna\'n\\ul. Uniwersytetu Pozna\'nskiego 4\\
61-614 Pozna\'n, Poland}
\email{piotr.mackowiak@amu.edu.pl}
\date{}
\title[Banach contraction principle and CH]{A converse of the Banach contraction principle for partial metric spaces and the continuum hypothesis}
\begin{document}
\begin{abstract}
A version of the Bessaga inverse of the Banach contraction principle for partial metric spaces is presented. Equivalence of that version and the continuum hypothesis is shown as well.
\end{abstract}
\maketitle
\section{Introduction}
The famous Bessaga converse of the Banach contraction principle (briefly: BCP) theorem  was published in 1959 \cite{B}. Loosely speaking, it states, for a map $T$ acting from a set $X$ to itself and possessing exactly one periodic point, that it is possible to make the map a contraction with respect to some special complete metric on $X$. This result by Bessaga gave rise to further investigations on the existence and properties of metric spaces that make a map like $T$ a contraction. There have appeared papers in which some topological conditions were imposed both on the map and the new metric \cite{HS, LJ, SL, PM}. Moreover, a few variants (or ways of presentation) of proofs of the Bessaga theorem have been published (see e.g. \cite[pp. 191--192]{KD}, \cite{JJ}, \cite[pp. 7--8]{VP}, \cite[pp. 525-526]{ES}). 
 
On the other hand, in 1992 the notion of a partial metric space was introduced by Matthews with the motivation, as stated in \cite{M}:
\begin{quotation} [...] to develop metric based tools for program verification in which the notion of the size of an object in a domain plays a pivotal role in quantifying the extent of its definedness.
\end{quotation}
In contrast to the case of a metric, in a partial metric space it may happen that the (partial) distance of an object to itself, the size of the object, is positive; an object is complete if its size is $0$. In the same work, Matthews has presented a version of BCP for partial metric spaces. Since the time of Matthews' work \cite{M} publication, there have emerged many works investigating purely mathematical side of partial metric spaces. For example, there have appeared variants of BCP for partial metric spaces \cite{ASS, BR, IPR, IR}, results concerning the issue of completeness of a partial metric space \cite{GL, HWZ, MM}, or, only recently, theorems relating sequential compactness and compactness of a partial metric space \cite{BMW, MM}.  

As we have mentioned, the paper \cite{M} contains a version of BCP for partial metric spaces. It seems natural to ask whether there is a Bessaga--type converse of the Banach contraction principle in partial metric setting. As long as one considers maps that have exactly one periodic point the answer is clear: since a metric space is also a partial metric space, an answer given in a metric setting also addresses the question in a partial metric setting. In this paper, basing on a special version of BCP formulated by  D.Ili\'{c}, V.Pavlovi\'{c} and V.Rako\v{c}evi\'{c} in \cite{IPR}, we show that a version of Bessaga's theorem on the inverse of BCP holds in a partial metric setting as well. Let us emphasize an essential fact that the version of BCP presented in \cite{IPR} does not exclude multiplicity of fixed points. To be more accurate, a map may have continuum of fixed points and still meet assumptions of BCP in partial metric setting \cite{IPR}. Interestingly, the inverse of BCP for partial metric spaces we propose is equivalent to the validity the continuum hypothesis (in short, CH).

The next section presents necessary notations and notions related to partial metric spaces. Section \ref{Results} contains main results of the paper. It includes a new version of proof of the Bessaga theorem, the statement and the proof of BCP inverse in partial metric setting. Finally, the equivalence of the latter theorem and CH is formulated and proved.

\section{Preliminaries}
We denote by $\bn$ the set of positive integers, $\bn_0:=\bn\cup\{0\}$. 
The set of non-negative real numbers is denoted by $\br_+:=\{x\in \br:\, x\geq 0\}$, $\aleph_0$ is the cardinality of $\bn$ and $\mathfrak{c}=2^{\aleph_0}$ is the cardinality of $\br$.

Let a set $X$ and a map $T:X\to X$ be given. For any point $x\in X$ we put $T^0(x):=x,T^1(x):=T(x),\ldots, T^n(x):=(\underbrace{T\circ\ldots\circ T}_{n\text{ times}})(x)$, $n\in \bn$,  $\, x\in A$. The orbit of a point $x\in X$ (under the map $T$) is the set $O^T(x):=\{x^T_n:\, n\in\bn_0\}\subset X$, where $x^T_n:=T^n(x),\, n\in \bn_0$; if there is no ambiguity we write $O(x)$ and $x_n$ instead of $O^T(x)$ and $x^T_n$, respectively. If a point $x\in X$ is a periodic point of $T$, then $O^T(x)$ is a finite set; otherwise it is infinite. By $p_{x'}(x)$ we denote the first $n\in \bn_0$ for which $x^T_n\in O(x')$; if there is no such $n$, then $p_{x'}(x):=+\infty$. It is clear that $p_{x'}(x)=p_{x'}(T(x))+1$ if $x\notin O^T(x')$. Moreover, $O^T(x)\subset O^T(x')$ for $x\in O^T(x')$.

Fix now a nonempty set $X$. A function $p:X\times X\to \br_+$ that satisfies, for any $x,z,y\in X$, the conditions
\begin{enumerate}[label={\textup{(\arabic*)}},ref=\textup{.(\arabic*)}]
\item $x=y\Leftrightarrow p(x,x)=p(x,y)=p(y,y)$;
\item $p(x,x)\leq p(y,x)$;
\item  $p(x,y)=p(y,x)$;
\item $p(x,y)\leq p(x,z)+p(z,y)-p(z,z)$.
\end{enumerate}
is called a partial metric on $X$. The pair $(X,p)$ is then called a partial metric space. If $p(x,x)=0$, $x\in X$, then the partial metric $p$ is a metric and the pair $(X,p)$ is a metric space. For a given partial metric space $(X,p)$ define (cf. \cite{IPR})
$$\rho_p:=\inf\{p(x,x):\, x\in X\},\quad X_p:=\{x\in X:\, p(x,x)=\rho_p\}.$$
Observe that it may happen that $X_p=\emptyset$.

Let us fix a partial metric space $(X,p)$. We say that a sequence $(x_n)_{n\in \mathbb{N}}$ of elements of $X$ converges to $x\in X$ if $\lim_{n\rightarrow \infty}p(x_n,x)=p(x,x)$. As usually, in this case the point $x$ is called a limit of the sequence $(x_n)_{n\in \mathbb{N}}$ and the sequence is said to be convergent. In general, partial metric spaces are not Hausdorff spaces and the limit may not be unique.
A sequence $(x_n)_{n\in \mathbb{N}}$ in $X$ properly converges to $x\in X$, if it converges to $x$ and $\lim_{n\to \infty}p(x_n,x_n)=p(x,x)$. It is clear that if $\lim_{n\to \infty}p(x_n,x)=p(x,x)$ and $\lim_{n\to \infty}p(x_n,x_n)=p(x,x)$, then $\lim_{m,n\to \infty}p(x_m,x_n)=p(x,x)$. There exist partial metric spaces such that some convergent series in these spaces do not converge properly. A sequence $(x_n)_{n\in \mathbb{N}}$ in $X$ is a Cauchy sequence if there exists the limit $\lim_{n,m\rightarrow\infty}p(x_n,x_m)\in \br_+.$ The partial metric space $(X,p)$ is complete if every Cauchy sequence $(x_n)_{n\in \mathbb{N}}$ of elements of $X$ is properly convergent. The partial metric space $(X,p)$ is $0$-complete if every Cauchy sequence $(x_n)_{n\in \mathbb{N}}$ elements of $X$ such that $\lim_{n,m\to\infty}p(x_n,x_m)=0$ properly converges to some element $x\in X$.

The above definitions related to partial metric spaces, many references and interrelationships can be found in \cite{CO}.


\section{Results}\label{Results}
We shall use the following lemmas first of which is obvious while the second is a small modification of a lemma from \cite{JJ}. 
\begin{lm}\label{lm:1a}
For all $a,b,c\in \br$, $\max\{a,b\}\leq \max\{a,c\}+\max\{c,b\}-c$.
\end{lm}
\begin{lm}\label{lm:1}
Let $X$ be any nonempty set and $\alpha \in (0,1)$ be an arbitrarily fixed number. If a map $T:X\to X$ possesses exactly one fixed point, say $\ov{x}\in X$, the following statements are equivalent:
\begin{enumerate} \item[\textup{(1)}] there exists a complete metric $d$ on $X$ such that
$d(T(x),T(y))\leq \alpha d(x,y),\,x,y\in X$, that is, $T$ is $\alpha$-contraction on the complete metric space $(X,d)$;
\item[\textup{(2)}] there exists a map $\varphi:X\to \br_+$ such that $\varphi^{-1}(0)=\{\ov{x}\}$ and $\varphi(T(x)) \leq \alpha \varphi(x)$, $x\in X$.
\end{enumerate}
\end{lm}
\begin{proof}
If $d$ is a complete metric on $X$ for which (1) is true, it suffices to put $\varphi(x):=d(x,\ov{x}), x\in X$. If $\varphi$ is as in assertion (2), then define $d(x,y):=\max\{\varphi(x),\varphi(y)\}$, $x\neq y,\, x,y\in X$, and $d(x,x):=0,\,x\in X$. By Lemma \ref{lm:1a} it is clear that $d$ is a metric on $X$. A sequence $(x_n)_{n\in \bn}$ in $(X,d)$ is Cauchy if and only if $(\varphi(x_n))_{n\in \bn}$ is either eventually constant or arbitrarily close to $0$ for large $n$. In the latter case, $\lim_{n\to \infty}d(x_n,\ov{x})=\lim_{n\to \infty}\max\{\varphi(x_n),\varphi(\ov{x})\}=\lim_{n\to \infty}\max\{\varphi(x_n),0\}=0=d(\ov{x},\ov{x})$. So, the metric space $(X,d)$ is complete. That $d(T(x),T(y))\leq \alpha d(x,y),\, x,y\in X$, is obvious.
\end{proof}
\subsection{Bessaga theorem in metric setting}
Here we present a proof of the Bessaga theorem. Our proof is novel to a limited extent only, however, the main reason to keep it here is that we use it in the later part of the paper.
\begin{thrm}[cf. \cite{B}]\label{thm:1} For a nonempty set $X$ and a map $T:X\to X$, let there be exactly one point $\ov{x}\in X$ for which $\ov{x}=T(\ov x)$. Suppose also that $x\neq T^n(x),\, n\in \bn$, $x\in X\backslash \{\ov{x}\}$. For each $\alpha\in (0,1)$ there exists a complete metric $d$ on $X$ such that $d(T(x),T(y))\leq \alpha d(x,y)$, $x,y\in X$.
\end{thrm}
\begin{proof} We shall prove that claim $(2)$ of Lemma \ref{lm:1} is satisfied. Recall that, for any $x\in X$, $x_n:=T^n(x)$, $n\in \bn_0$, and $O(x):=\{x_n:\, n\in \bn_0\}$. Observe that the only periodic point of $T$ is the unique fixed point $\ov{x}$. For $x\in X$, denote $V_x:=\bigcup\{O(y):\, y\in X\text{ and }O(y)\cap O(x)\neq\emptyset\}$. By this definition, for any $x, x'\in X$, we have $x\in V_x$, and either $V_{ x}\cap V_{ x'}=\emptyset$ or $V_{x}=V_{x'}$. Indeed, suppose that $z\in V_x\cap V_{x'}$. There are $y\in V_x, y'\in V_{x'}$ with $z\in O(y)$ and $z\in O(y')$. But $O(y)\cap O(x)\neq \emptyset$ and $O(y)\cap O(x)\neq \emptyset$. Therefore, $O(z)\cap O(x)\neq \emptyset$ and $O(z)\cap O(x')\neq \emptyset$. Thus, there exist  $m_z,n_z,n'_z\in \bn_0$ such that $z_{m_z+q}=x_{n_z+q}=x'_{n'_z+q}$ for all $q\in \bn_0$. Now, let $a\in V_x$. It follows that $a_{m_a+r}=x_{n_a+r}$ for some $m_a,n_a\in \bn$ and all $r\in \bn_0$. It suffices to choose $q,r\in \bn_0$ such that $m_z+q=m_a+r$ (which is always possible) to get $a_{m_a+r}=x_{n_z+q}=x'_{n'_z+q}$ and, in consequence, $a\in V_{x'}$. This shows that $V_x=V_{x'}$, if $V_x\cap V_{x'}\neq \emptyset$. 

By the above considerations, $X=\bigcup_{\wh x\in W}V_{\wh x}$, for a set $W\subset X$ such that $V_{\wh x}\cap V_{\wh x'}=\emptyset$ if $\wh x\neq\wh x'$, $\wh x,\wh x'\in W$. The existence of the set $W$ is a consequence of the axiom of choice: the relation $\sim$ defined by $x\sim y\Leftrightarrow V_x=V_y$ is an equivalence relation on $X$ and to construct $W$ we just pick an element from each equivalence class of $\sim$ and include it into $W$. Without loss of generality we also assume that $\ov{x} \in W$.

Let us define $\varphi_{\ov{x}}:V_{\ov x}\to \br_+$ by $\varphi_{\ov x}(x):=\alpha^{-p_{\ov{x}}(x)}$, $x\in X\backslash\{\ov x\}$, and $\varphi_{\ov x}(\ov x):=0$. Since $O(\ov x)=\{\ov x\}$ for $x\in V_{\ov x}\backslash\{\ov x\}$, we have $p_{\ov{x}}(T(x))=p_{\ov{x}}(x)-1$. Consequently, $\varphi_{\ov{x}}(T(x))\leq \alpha^{-(p_{\ov{x}}(x)-1)}= \alpha \varphi_{\ov{x}}(x), x\in V_{\ov{x}}\backslash\{\ov x\}$, with strict inequality only if $T(x)=\ov{x}$. Obviously, $\varphi_{\ov{x}}(T(\ov{x}))=0\leq \alpha \varphi_{\ov{x}}(\ov x)=0$.

For $\wh x\in W\backslash \{\ov x\}$, define $\varphi_{\wh x}:V_{\wh x}\to \br_+$ by 
\begin{equation}\label{eqn:0}\varphi_{\wh x}(x):=\alpha^{-p_{\wh x}(x)+p_{w(x)}(\wh x)},\end{equation}
where $w(x):=T^{p_{\wh x}(x)}(x)$, $x\in V_{\wh x}$; $w(x)$ is simply the first term of the sequence $x_0,x_1,\ldots$ that belongs to $O(\wh x)$. Moreover, $p_{\wh x}(\wh x_n)=0$ (by definition of $p_{\wh x}(\cdot)$), $p_{\wh x_n}(\wh x)=n$ (since $\wh x $ is not periodic) which results in $w(\wh x_n)=\wh x_n$, $n\in \bn_0$. 

Formula \pref{eqn:0} and equality $\wh x_n=w(\wh x_n)$ imply $\varphi_{\wh x}(\wh x_n)=\alpha^{-p_{\wh x}(\wh x_n)+p_{w(\wh x_n)}(\wh x)}=\alpha^{p_{w(\wh x_n)}(\wh x)}=\alpha^n$, $n\in \bn_0$. Hence, $\varphi_{\wh x}(T(\wh x_n))=\varphi_{\wh x}(\wh x_{n+1})=\alpha^{n+1}=\alpha\alpha^n=\alpha \varphi_{\wh x}(\wh x_n)$ as $n\in \bn_0$. For $x\in V_{\wh x}\backslash O(\wh x)$, we obtain $$\varphi_{\wh x}(T(x))=\alpha^{-p_{\wh x}(T(x))+p_{w(T(x))}(\wh x)}\\=\alpha^{-(p_{\wh x}(x)-1)+p_{w(T(x))}(\wh x)}=\alpha^{-(p_{\wh x}(x)-1)+p_{w(x)}(\wh x)}=\alpha \varphi_{\wh x}(x),$$ where we used the equality $w(T(x))=w(x)$, $x\in V_{\wh x}\backslash O(\wh x)$.
Now, let $\varphi:X\to\br_+$ be defined by $\varphi(x):=\varphi_{\wh x}(x),$ where $\wh x\in W$ is the unique element of $W$ with $x\in V_{\wh x}$. It follows that $\varphi$ satisfies the conditions of Lemma \ref{lm:1}.(2).
\end{proof}
\subsection{Bessaga theorem in partial metric setting}
In \cite{IPR} there is stated and proved a theorem which we find as a natural counterpart of the Banach contraction principle for partial metric spaces.
\begin{thrm}[cf. Theorem 3.1 in \cite{IPR}]\label{thm:2}
Let $(X,p)$ be a complete partial metric space and suppose that $T:X\to X$ is a map satisfying the following condition  
$$p(T(x),T(y))\leq \max\{\alpha p(x,y),p(x,x),p(y,y)\},$$
$x,y\in X$, where $\alpha \in (0,1)$ is fixed. Then there exists a unique point $\ov{x}\in X_p$ for which $T(\ov{x})=\ov{x}$. Moreover, $\lim_{n\to\infty}p(T^n(x),\ov{x})=p(\ov{x},\ov{x})=\lim_{m,n\to\infty}p(T^n(x),T^m(x))$ for any $x \in X_p$.
\end{thrm}
By slightly modifying the original proof of Theorem \ref{thm:2} one can prove the next theorem.  
\begin{thrm}\label{thm:3}
Let $(X,p)$ be a $0$--complete partial metric space with $\rho_p=0$. Suppose that $T:X\to X$ is a map satisfying the following condition
\begin{equation}\label{eqn:1}p(T(x),T(y))\leq \max\{\alpha p(x,y),p(x,x),p(y,y)\},\end{equation}
$x,y\in X$, where $\alpha \in (0,1)$ is fixed. Then there exists a unique point $\ov{x}\in X_p$ for which $T(\ov{x})=\ov{x}$. Moreover, $0=p(\ov{x},\ov{x})=\lim_{n\to\infty}p(T^n(x),\ov{x})=\lim_{m,n\to\infty}p(T^n(x),T^m(x))$ for any $x \in X_p$.
\end{thrm}
Notice that if a map $T$ satisfies condition \pref{eqn:1}, then for any $a\in \br $ there exists at most one fixed point $x_a$ of $T$ for which $p(x_a,x_a)=a$. Therefore, any such a map has at most continuum of fixed points, that is, $\textup{card}\{x\in X:\, x=T(x)\}\leq \mathfrak{c}$.

It turns out that a version of the Bessaga theorem holds true in partial metric spaces.
\begin{thrm}\label{thm:4}
Let $X$ be a nonempty set and $T:X\to X$ be a map such that the nonempty set $A:=\{x\in X:\, x=T(x)\}$ of its fixed points is of cardinality at most continuum, that is, $\textup{card}(A)\leq \mathfrak{c}$. Suppose also that $x\neq T^n(x),\, n\in \bn, x\in X\backslash A$. For each $\alpha\in (0,1)$ there exists a partial metric $p$ on $X$ such that $(X,p)$ is a $0$-complete partial metric space with $\rho_p=0$ and it holds
\begin{equation}\label{eqn:2}p(T(x),T(y))\leq \max\{\alpha p(x,y),p(x,x),p(y,y)\},\end{equation}
$x,y\in X$. If CH holds, that is, $\textup{card}(A)=\mathfrak{c}$ if $\aleph_0<\textup{card}(A)\leq \mathfrak{c}$, then the space $(X,p)$ is a complete partial metric space.
\end{thrm}
\begin{proof}
We shall prove the claim in a way similar to the method presented in the proof of Theorem \ref{thm:1}. However, this time we do not have Lemma \ref{lm:1} at hand. As previously, for $x\in X$, denote $V_x:=\bigcup\{O(y):\, y\in X\text{ and }O(y)\cap O(x)\neq\emptyset\}$. Choose any $\ov{x}\in A$. Define $X_1:=\bigcup_{x\in A\backslash\{\ov x\}}V_x$ and $X_0:=X\backslash X_1$. Observe that $T(X_i)\subset X_i$, $i=0,1$, and $X_0\cap X_1=\emptyset$. Moreover, there is exactly one fixed point of $T$ that belongs to $X_0$, namely $\ov{x}$, and there are no other periodic points of $T$ in $X_0$. Therefore, $X_0=V_{\ov x}\cup \bigcup_{\wh x\in W}V_{\wh x}$, for a set $W\subset X_0\backslash\{\ov x\}$, and $V_x\cap V_y=\emptyset$ for $x\neq y$, $x,y\in A\cup W$ (see the proof of Theorem \ref{thm:1}). Let $h:A\backslash\{\ov{x}\}\to (\alpha,1)$ be an injection -- it exists by the assumption that $\textup{card}(A)\leq \mathfrak{c}$. If $\textup{card}(A)\leq \aleph_0$, then we assume that either $h(A\backslash\{\ov{x}\})$ is discrete or it has exactly one accumulation point in $(\alpha,1)$. Under the continuum hypothesis, if $\textup{card}(A)=\mathfrak{c}$ we may assume that $h(A\backslash\{\ov x\})$ is a closed subinterval of $ (\alpha+\varepsilon, 1-\varepsilon)$ for some $\varepsilon>0$. 

We shall now define a function $\varphi:X\to \br_+$ with help of which we will construct a partial metric on $X$. For any $x\in X$, by $\wh x$ we denote the unique $y\in A\cup W$ with $x\in V_y$. Hence, $\wh x=x$, $x\in A$. For $x\in X$ we define a map $\varphi:X\to \br_+$ by
$$\varphi(x):=\left\{\begin{array}{ll}
0&\text{if } x=\ov{x},\\
\alpha^{-p_{\ov{x}}(x)}&\text{if } x\in V_{\ov x}\backslash \{\ov x\},\\
\alpha^{-p_{\wh x}(x)+p_{w(x)}(\wh x)}&\text{if } x\in X_0\backslash V_{\ov x},\\
\alpha^{-p_{\wh{x}}(x)}h(\wh x)&\text{if } x\in X_1.
\end{array}\right.$$
It follows that $\varphi(T(x))\leq \alpha \varphi(x)$, $x\in X\backslash A$, with strict inequality only if $T(x)=\ov{x}$ (see the proof of Theorem \ref{thm:1} for details). Moreover, $\varphi(x)\geq \alpha$, $x\in X_1$.

Let us now define a function $p:X\times X\to \br_+$ by 
$$p(x,y):=\left\{\begin{array}{ll}\max\{\varphi(x),\varphi(y)\}&\text{if } x\neq y, x,y\in X\backslash A,\\
\max\{\frac{1}{2}\varphi(x),\frac{1}{2}\varphi(y)\}&\text{if }  x\neq y, x,y\in A,\\
\max\{\frac{1}{2}\varphi(x),\varphi(y)\}&\text{if }  x\neq y, x\in A, y\in X\backslash A,\\
\max\{\varphi(x),\frac{1}{2}\varphi(y)\}&\text{if }  x\neq y, x\in X\backslash A, y\in A,\\
0&\text{if }  x=y, x\in X_0,\\
\frac{1}{2}\varphi(x)&\text{if }  x=y, x\in X_1.
\end{array}\right.$$
The function $p$ is a partial metric on $X$. It is obvious that $p$ is nonnegative and symmetric. Observe that $p(x,x)\leq \frac{1}{2}\varphi(x)\leq \max\{\frac{1}{2}\varphi(x),\frac{1}{2}\varphi(y)\}\leq p(x,y)$, provided $x\neq y$, $x,y\in X$. Therefore, $p(x,x)\leq p(x,y),\, x,y\in X$. To prove the triangle inequality for partial metrics assume that $x,y,z\in X$. If $x=y$, then $p(x,y)=p(x,x)\leq p(x,z)\leq p(x,z)+p(z,y)-p(z,z)$. If $x=z$, then $p(x,y)=p(z,y)\leq p(x,z)+p(z,y)-p(z,z)$; the case $z=y$ is now obvious. So, let $x\neq z\neq y\neq x$. There are $a, b,c\in\{\frac{1}{2}, 1\}$ for which $p(x,y)=\max\{a\varphi(x),b\varphi(y)\}$, $p(x,z)=\max\{a\varphi(x),c\varphi(z)\}$ and $p(z,y)=\max\{c\varphi(z),b\varphi(y)\}$. Thus, by Lemma \ref{lm:1a} and nonnegativity of $\varphi$, $p(x,y)=\max\{a\varphi(x),b\varphi(y)\}\leq \max\{a\varphi(x),c\varphi(z)\}+\max\{c\varphi(z),b\varphi(y)\}-c\varphi(z)\leq p(x,z)+p(z,y)-\frac{1}{2}\varphi(z)=p(x,z)+p(z,y)-p(z,z)$. 

We now show that $p(x,y)=p(x,x)=p(y,y)$ implies $x=y$. If $p(x,x)=0$, then $x,y\in X_0$, but $p$ restricted to $X_0\times X_0$ is a metric on $X_0$ (see the proof of Theorem \ref{thm:1}), so $x=y$. Suppose that $p(x,x)>0$. Hence, $x,y\in X_1$ and $\varphi(x)>0, \varphi(y)>0$ which entails that $p(x,y)=\max\{\frac{1}{2}\varphi(x),\frac{1}{2}\varphi(y)\}>0$ and, consequently, $x,y\in A\backslash \{\ov{x}\}$. Now, if $x\neq y$, then $h(x)\neq h(y)$ and, equivalently, $\varphi(x)\neq\varphi(y)$ which is impossible. Therefore, $x=y$ and we have just shown that the pair $(X,p)$ is a partial metric space. 

To prove $0$-completeness notice that $\lim_{n,m\to \infty}p(x_n,x_m)=0$, for a sequence $(x_n)_{n\in\bn}$ in $X$, implies either the sequence is eventually constant with $x_n\in X_0$ for large $n\in \bn$ or $\lim_{n\to \infty}\varphi(x_n)=0$ which implies $x_n\in X_0$ and $0=p(\ov{x},\ov{x})=\varphi(\ov{x})=\lim_{n\to\infty} p(x_n,\ov{ x})=\lim_{n\to \infty}p(x_n,x_n)$. Whence, the partial metric space $(X,p)$ is $0$-complete. Suppose now that CH is true and $\lim_{n,m\to \infty}p(x_n,x_m)=a>0$ for a sequence $(x_n)_{n\in\bn}$ in $X$. Then we may assume that $x_n\in X_1,\, n\in \bn$. This implies that $\lim_{n\to \infty} p(x_n,x_n)=\lim_{n\to \infty} \frac{1}{2}\varphi(x_n)=a$. Since $a>0$, we see that, for large $m,n\in \bn$, $\max\{\frac{1}{2}\varphi(x_n), \varphi(x_m)\}>\frac{3}{2}a$, $\max\{\varphi(x_n),\frac{1}{2}\varphi(x_m)\}>\frac{3}{2}a$ and, the more, $\max\{\varphi(x_n),\varphi(x_m)\}>\frac{3}{2}a$. This, in view of the definition of the partial metric $p$, proves that either the sequence $(x_n)_{n\in \bn}$ is eventually constant (because $\lim_{n\to\infty} p(x_n,x_n)=\lim_{n\to\infty}\frac{1}{2}\varphi(x_n)=a>0$) or $x_n\in A\backslash\{\ov{x}\}$ for large $n$. We may assume that $x_n\in A\backslash \{\ov x\}$, $n\in\bn$. 
It suffices to consider the latter case to prove completeness, because the former one is obvious. Notice that $\wh x=x$ and $p_{\wh{x}}(x)=0$ for $x\in A$. Hence, $p(x_n,x_m)=\max\{\frac{1}{2}\varphi(x_n),\frac{1}{2}\varphi(x_m)\}=\max\{\frac{1}{2}h(x_n),\frac{1}{2}h(x_m)\}$ and, since the sequence $(x_n)_{n\in \bn}$ is Cauchy, $a=\lim_{n\to \infty}\frac{1}{2}h(x_n)=\frac{1}{2}\wh h$ for some $\wh h\in h(A\backslash\{\ov{x}\})$, due to closedness of the image $h(A\backslash\{\ov{x}\})$ as a subset of $[0,1]$. Let $\wh h= h(x)$ for some (unique) $x\in A\backslash \{\ov x\}$. We have $\lim_{n\to\infty} p(x_n,x)=\lim_{n\to\infty} \max\{\frac{1}{2}h(x_n),\frac{1}{2}h(x)\}=\frac{1}{2}h(x)=p(x,x)=a$, which, taking into account that $\lim_{n\to \infty}p(x_n,x_n)=a$, shows that $(X,p)$ is a complete partial metric space.
To finish the proof it suffices to verify condition \pref{eqn:2}. Let $x,y\in X\backslash A$. Then, for some $a,b\in\{0,\frac{1}{2},1\}$, $p(T(x),T(y))=\max\{a\varphi(T(x)),b\varphi(T(y))\}\leq \alpha \max\{a\varphi(x),b\varphi(y)\}\leq \alpha p(x,y)$. If $x\in A$ and $y\in X\backslash A$, then $\frac{1}{2}\varphi(T(x))=\frac{1}{2}\varphi(x)=p(x,x)$ and $\varphi(T(y))\leq \alpha \varphi(y)\leq \alpha p(x,y)$, from which we get $p(T(x),T(y))\leq \max\{\frac{1}{2}\varphi(T(x)),\varphi(T(y))\}\leq \max\{\alpha p(x,y), p(x,x)\}$ (the case $y\in A$ and $x\in X\backslash A$ is symmetric). For $x,y\in A$, $p(T(x),T(y))=\max\{\frac{1}{2}\varphi(T(x)),\frac{1}{2}\varphi(T(y))\}=\max\{\frac{1}{2}\varphi(x),\frac{1}{2}\varphi(y)\}=\max\{p(x,x),p(y,y)\}$. The claim follows.
\end{proof}

\subsection{Bessaga theorem for partial metric spaces and CH}
If CH is true, the partial metric space constructed in the proof of Theorem \ref{thm:4} is complete. However, if the set of fixed points, that is, $A$ is uncountable and we do not know whether it is of cardinality continuum, then it may happen that $h(A)$ is not bijective with any closed subset of $\br$. Recall that according to the Cantor-Bendixson theorem uncountable closed subsets of $\br$ contain perfect subsets (see \cite{K} or \cite{KK} for set--theoretical or topological definitions and theorems used in this subsection). It is well-known that perfect subsets of $\br$ are of cardinality $\mathfrak{c}$. In view of Theorem \ref{thm:4}, it turns out that CH is equivalent to the statement: for any fixed $\alpha\in (0,1)$, if $A\neq \emptyset$ is a set with $\textup{card}(A)\leq \mathfrak{c}$, then there exists a partial metric $p$ on $A$ such that the partial metric space $(A,p)$ is complete and the identity mapping $Id(x):=x$, $x\in A$, satisfies the condition $p(Id(x), Id(y))\leq \max\{\alpha p(x,y),p(x,x),p(y,y)\}$, $x,y\in A$. Observe that the inequality $p(Id(x), Id(y))\leq \max\{\alpha p(x,y),p(x,x),p(y,y)\}$, $x,y\in A$, is equivalent to $p(x,y)\leq \max\{p(x,x),p(y,y)\}$, $x,y\in A$. We summarize the above discussion in the last theorem in this paper.

\begin{thrm}\label{thm:equiv}
The following statements are equivalent.
\begin{enumerate}[label={\textup{(\arabic*)}},ref=\textup{(\arabic*)}] 
\item\label{thm:equiv:1} CH holds, that is, for any set $A$, if $\aleph_0<\textup{card}(A)\leq \mathfrak{c}$, then $\textup{card}(A)=\mathfrak{c}$.
\item\label{thm:equiv:1a} Let $X$ be a nonempty set, $T:X\to X$ a map such that $0<\textup{card}\{x\in X:\, x=T(x)\}\leq \mathfrak{c}$ and $x\neq T^n(x),\, n\in \bn, x\in X\backslash \{x\in X:\, x=T(x)\}$. Then, for each $\alpha\in (0,1)$, there exists a partial metric $p$ on $X$ such that $(X,p)$ is a complete partial metric space with $\rho_p=0$ and
$$p(T(x),T(y))\leq \max\{\alpha p(x,y),p(x,x),p(y,y)\},\quad x,y\in X.$$
\item\label{thm:equiv:2} For each set $A$ whose cardinality does not exceed $\mathfrak{c}$ there exists a partial metric $p$ on $A$ such that the partial metric space $(A,p)$ is complete, $\rho_p=0$, and $p(x,y)\leq \max\{p(x,x),p(y,y)\}$, $x,y\in A.$
\end{enumerate}
\end{thrm}
\begin{proof} By Theorem \ref{thm:4}, \ref{thm:equiv:1} implies \ref{thm:equiv:1a}. Setting $T(x):=Id(x), \, x\in X:=A$, we get that \ref{thm:equiv:1a} entails \ref{thm:equiv:2} (see the discussion preceding Theorem \ref{thm:equiv}). 

We shall show that \ref{thm:equiv:1} follows from \ref{thm:equiv:2}. Suppose that $(A,p)$ is a complete partial metric space such that $p(x,y)\leq \max\{p(x,x),p(y,y)\}$, $x,y\in A$. It suffices to show $\aleph_0< \textup{card}(A)\leq \mathfrak{c}$ implies $\textup{card}(A)=\mathfrak{c}$. Let $f(x):=p(x,x),\, x\in A$. Observe that $f(x)\neq f(y)$ for $x\neq y$ - this follows from definition of a partial metric. Hence, $f$ is injective and it maps $A$ onto $f(A)\subset \br$. We shall show that $f(A)$ is a closed subset of $\br$. Let $(y_n)_{n\in \bn}$ be a sequence in $f(A)$ that converges to some $y\in \br$. For each $n\in \bn$ there is $x_n\in A$ with $f(x_n)=y_n$. Observe that $y_n=f(x_n)=p(x_n,x_n)\leq p(x_n,x_m)\leq \max\{p(x_n,x_n),p(x_m,x_m)\}=\max\{f(x_n),f(x_m)\}=\max\{y_n,y_m\}$. Hence, the sequence $(x_n)_{n\in \bn}$ is a Cauchy sequence in $(A,p)$ with $\lim_{m,n\to\infty}p(x_n,x_m)=y$. By completeness of $(A,p)$, there is $x\in A$ such that $\lim_{n\to\infty}p(x_n,x)=p(x,x)=y$. From this it follows that $y\in f(A)$ and $f(A)$ is a closed uncountable subset of $\br$. By the Cantor-Bendixson theorem, $\mathfrak{c}\leq\textup{card}(f(A))\leq \mathfrak{c}$. Equality $\textup{card}(f(A))=\mathfrak{c}$ follows from the Cantor-Bernstein theorem. But $A$ and $f(A)$ are of the same cardinality, so $\textup{card}(A)=\mathfrak{c}$.
\end{proof}



\begin{thebibliography}{99}
\bibitem{ASS} I.Altun, H.Simsek, F.Sola, \emph{Generalized contractions on partial metric spaces}, Topology Appl. 157, 2010, 2778--2785.


\bibitem{B} C.Bessaga, \emph{On the converse of the Banach fixed-point principle}, Colloq. Math. 7, 1959, 41--43.

\bibitem{BMW} D.Bugajewski, P.Ma\'ckowiak, R.Wang, \emph{On compactness and fixed point theorems in partial metric spaces}, Fixed Point Theory, 23, 2023, 163--178.

\bibitem{BR} D.Bugajewski, R.Wang, \emph{On the topology of partial metric spaces}, Math. Slovaca, 70, 2020, 135--146.

\bibitem{CO} S.Cobza\c{s}, \emph{Fixed points and completeness in metric and in generalized metric spaces}, J. Math. Sci. 250, 2020, 475--535.
\bibitem{KD} K.Deimling, \emph{Nonlinear functional analysis}, Berlin, Springer, 1985.

\bibitem{GL} X.Ge, S.Lin, \emph{Completions of partial metric spaces}, Topology Appl. 182, 2015, 16--23.

\bibitem{HWZ} S.Han, J.Wu, D.Zhang, \emph{Properties and principles on partial metric spaces}, Topology Appl., 230, 2017, 77--98.

\bibitem{HS} P.Hitzler, A.K.Seda, \emph{A ``converse'' of the Banach contraction mapping theorem}, Journal of Electrical Engineering 52, 2001, 3--6.

\bibitem{IPR} D.Ili\'{c}, V.Pavlovi\'{c} and V.Rako\v{c}evi\'{c}, \emph{Some new extensions of Banach's contraction principle to partial metric space}, Appl. Math. Lett. 24, 2011, 1326--1330.

\bibitem{JJ} J.Jachymski, \emph{A short proof of the converse to the contraction principle and some related results}, Topol. Methods Nonlinear Anal. 15, 2000, 179--186.
\bibitem{LJ} L.Janos, \emph{A converse of Banach's contraction theorem}, Proc. Amer. Math. Soc. 16, 1967, 287--289.
\bibitem{K} A.S.Kechris, \emph{Classical descriptive set theory}, Springer-Verlag, New York, 1985.
\bibitem{KK} K.Kuratowski, \emph{Introduction to set theory and topology}, Pergamon Press, Oxford, 1961.

\bibitem{SL} S.Leader, \emph{A topological characterization of Banach contractions}, Pacific J. Math. 69, 1977, 461--466. 

\bibitem{M} S.G.Matthews, \emph{Partial metric topology}, in: Papers on general topology and applications, Flushing, NY, 1992, in: Ann. New York Acad. Sci., vol.728, New York Acad.Sci., New York, 1994, 183--197.
\bibitem{PM} P.R.Meyers, \emph{A converse to Banach's contraction theorem}, J. Res. Natl. Bur. Stand. 71B, 1967, 73--76.

\bibitem{MM} V.Mykhaylyuk, V.Myronyk, \emph{Compactness and completeness in partial metric spaces}, Topology Appl. 270, 2020,  106925.

\bibitem{VP} V.Pata, \emph{Fixed point theorems and applications}, Springer, 2019

\bibitem{IR} I.A.Rus, \emph{Fixed point theory in partial metric spaces}, An. Univ. Vest Timi. s. Ser. Mat.-Inform. 46, 2008, 149--160.

\bibitem{ES} E.Schechter, \emph{Handbook of analysis and its foundations}, Academic Press, London, 1996.

\end{thebibliography}
\end{document}